\documentclass[12pt]{article}

\usepackage{graphicx}
\usepackage[usenames]{color}

\usepackage{amsfonts}
\usepackage{amsmath}
\usepackage{amssymb}
\usepackage{graphicx}
\usepackage[usenames]{color}

\newtheorem{theorem}{Theorem}[section]

\newtheorem{corollary}[theorem]{Corollary}
\newtheorem{lemma}[theorem]{Lemma}
\newtheorem{remark}[theorem]{Remark}
\newtheorem{definition}[theorem]{Definition}

\def\hf{{\hfill\rule{1,5mm}{1,5mm}}}
\def\1{^{-1}}
\def\ff{\forall }

\def\oo{{\omega}}

\def\dd{\displaystyle}
\def\ooo{{\Omega}}
\def\<{\left<}
\def\>{\right>}
\def\({\left(}
\def\){\right)}
\def\ff{\forall }
\def\D{\Delta }
\def\9{{\infty}}
\def\barr{\begin{array}}
\def\earr{\end{array}}
\def\ov{\overline}
\def\wt{\widetilde}
\def\vp{{\varepsilon}}

\def\pp{{\partial}}
\def\dd{\displaystyle}
\def\vf{{\varphi}}
\def\lbb{{\lambda}}
\def\g{{\gamma}}
\def\a{{\alpha}}
\def\b{{\beta}}
\def\pas{\mathbb{P}\mbox{-a.s.}}
\def\de{{\delta}}
\def\cala{{\mathcal{A}}}
\def\calb{{\mathcal{B}}}
\def\cald{{\mathcal{D}}}
\def\calf{{\mathcal{F}}}

\def\calo{\mathcal{O}}
\def\calf{\mathcal{F}}
\def\calh{\mathcal{H}}

\def\calx{\mathcal{X}}

\def\calz{{\mathcal{Z}}}

\def\vsp{\vspace*{1,5mm}\\ }
\def\vspp{\vspace*{1,5mm}\\ }
\def\n{\noindent }

\def\D{{\Delta}}
\def\rr{{\mathbb{R}}}

\def\E{{\mathbb{E}}}
\def\nn{{\mathbb{N}}}
\def\dd{\displaystyle}

\def\3{\subset }

\def\sk{\smallskip }

\def\bk{\bigskip }

\begin{document}

\title{\bf Variational solutions to nonlinear stochastic differential equations in~Hilbert spaces}

\author{Viorel Barbu\thanks{Octav Mayer Institute of Mathematics of  Romanian Academy),     Ia\c si, Romania.  Email: vbarbu41@gmail.com}\and Michael R\"ockner\thanks{Fakult\"at f\"ur Mathematik, Universit\"at Bielefeld,  D-33501 Bielefeld, Germany.  Email: roeckner@math.uni-bielefeld.de}}
\date{}

\maketitle

\begin{abstract} One introduces a new variational concept of solution for the stochastic dif\-fe\-ren\-tial equation $dX+A(t)X\,dt+\lbb X\,dt=X\,dW,$ $t\in(0,T)$; $X(0)=x$ in a real Hilbert space where $A(t)=\pp\vf(t),$ $t\in(0,T)$, is a ma\-xi\-mal monotone subpotential operator in $H$ while $W$ is a   Wiener process in $H$ on a probability space $\{\ooo,\calf,\mathbb{P}\}$. In this new context, the solution $X=X(t,x)$ exists for each $x\in H$, is unique, and depends continuously on $x$. This functional scheme applies to a general class of stochastic PDE not covered by the classical variational existence theory (\cite{8}, \cite{13-b}, \cite{9}) and, in particular, to stochastic variational ine\-qua\-li\-ties and parabolic stochastic equations with general   monotone nonlinearities with low or superfast growth to $+\9$.\\
{\bf  Keywords:} Brownian motion, maximal monotone operator, subdifferential, random differential equation, minimizastion problem.\\
{\bf Mathematics Subject Classification (2010):} Primary 60H15; Secondary 47H05, 47J05.
\end{abstract}

\section{Introduction}\label{s1}
Here, for $\lbb\in(0,\9)$, we consider   the stochastic differential equation
\begin{equation}\label{e1.1}
\barr{ll}
dX(t)+A(t)X(t)dt+\lbb X(t)dt\ni X(t)dW_t,\ t\in(0,T),\vspp
X(0)=x\in H,\earr
\end{equation}in a real Hilbert spaced $H$ whose elements are generalized functions on a bounded domain $\calo\subset\rr^d$ with a smooth boundary $\pp\calo$. In examples, we have in mind that $H$ is e.g.  $L^2(\calo)$ or $H^1_0(\calo)$, $H^1(\calo)$, $H^{-1}(\calo)$.

The norm of $H$ is denoted by $|\cdot|_H$, its  scalar product by $(\cdot,\cdot)$ and its Borel $\sigma$-algebra by $\calb(H)$.

$W$ is a Wiener process of the form
\begin{equation}\label{e1.2}
W(t,\xi)=\sum^\9_{j=1}\mu_je_j(\xi)\beta_j(t),\ \xi\in\calo,\ t\ge0,\end{equation}
where $\{\beta_j\}^\9_{j=1}$ is an independent system of real $(\calf_t)$-Brownian motions on a probability space $\{\ooo,\calf,\mathbb{P}\}$ with natural filtration $(\calf_t)_{t\ge0}$ and $\{e_j\}$ is an orthonormal basis in $H$ such that both  $c_j$ and $e^2_j$, $j\in\nn$, are multipliers in $H$, while $\mu_j\in\rr,$ $j=1,2,...$, satisfy  \eqref{e1.9} below.

As regards the nonlinear (multivalued) operator $A=A(t,\oo):H\to H$, the following hypotheses  will be assumed below.
\begin{itemize}\item[(i)] {\it Let $\vf:[0,T]\times H\times\ooo\to \ov\rr=]-\9,+\9]$ be  convex lower semicontinuous in $y\in H$ and progressively measurable, i.e., for each $t\in[0,T]$ the  function $\vf$ restricted to $[0,t]\times H\times\ooo$ is $\calb([0,t])\otimes\calb(H)\otimes\calf_t$ measurable, and let
\begin{equation}\label{e1.3}
A(t,\oo)=\pp\vf(t,\oo),\ \ff(t,\oo)\in[0,T]\times\ooo.\end{equation}
 In particular, $y\to A(t,\oo,y)$ is maximal monotone in $H\times H$ for all $(t,\oo)\in[0,T]\times\ooo$. Furthermore, $\vf$ is such that there exists $\alpha\in L^2([0,T]\times\ooo;H)$ and $\beta\in L^2([0,T]\times\ooo)$ such that

 {$\vf(t,y,\oo)\ge(\alpha(t,\oo),y)-\beta(t,\oo)\mbox{ for }dt\otimes\mathbb{P}-a.e.,\ (t,\oo)\in[0,T]\times\ooo.$}}

\item[(ii)] {\it $e^{\pm W(t)}$ is a  multiplier in $H$ such that there is an $(\calf_t)_{t\ge0}$-adapted $\rr_+$-valued process $Z(t),$ $t\in[0,T]$, with
\begin{eqnarray}
&\barr{rcll}
\sup\limits_{t\in[0,T]}|Z(t)| &<&\9,&\pas,\vsp
|e^{\pm W(t)}y|_H&\le&Z(t)|y|_H,&\ff t\in[0,T],\ y\in H,\earr\label{e1.4}\\[0,2mm]
&t\to e^{\pm W(t)}\in L(H,H)\mbox{ is continuous.}\label{e1.5} \end{eqnarray}}
\end{itemize}

 Recall that a multivalued mapping $A:D(A)\subset H\to H$ is said to be maximal monotone if it is monotone, that is, for $u_1,u_2\in D(A)$,
$$(z_1-z_2,u_1-u_2)\ge0,\ \ff z_i\in Au_i,\ i=1,2,$$and the range $R(\lbb I+A)$ is all of $H$ for each $\lbb>0$.

If $\vf:H\to \ov \rr$ is a convex, lower semicontinuous function, then its sub\-dif\-fe\-rential $\pp\vf:H\to H$
\begin{equation}\label{e1.6}
\pp\vf(u)=\{v\in H;\ \vf(u)\le\vf(\bar u)+(v,u-\bar u),\ \ff\bar u\in H\}\end{equation}is maximal monotone (see, e.g., \cite{1}).

The conjugate $\vf^*$ of $H$ defined by
\begin{equation}\label{e1.7}
\vf^*(v)=\sup\{(u,v)-\vf(u);\ u\in H\}
\end{equation}satisfies
\begin{equation}\label{e1.8}
\barr{rcll}
\vf(u)+\vf^*(v)&\ge&(u,v),&\ff u,v\in H,\vspp
\vf(u)+\vf^*(v)&=&(u,v),&\mbox{iff }v\in\pp\vf(u).\earr
\end{equation}

As regards the basis $\{e_j\}$ arising in the definition of the Wiener process~$W$, we assume also that, for the multipliers $e^2_j$, we have
\begin{itemize}\item[(iii)] {\it For $\g_j=\max\{\sup\{|ue_j|_H;\ |u|_H=1\},(\sup\{ u e^2_j|_H;\,|u|_H=1\})^{\frac12},1\}$, we assume
\begin{equation}\label{e1.9}
\nu=\sum^\9_{j=1}\mu^2_j\g^2_j<\9,
\end{equation}and that $\lbb>\nu$. }\end{itemize}

Clearly, then
\begin{equation}\label{e1.10}
\mu=\frac12\sum^\9_{j=1}\mu^2_je^2_j
\end{equation}is a  multiplier in $H$.

It should be noted that the condition $\lbb>\nu$ in \eqref{e1.9} is made only for   convenience. In fact, by the substitution $X\to\exp(-\lbb t)X$ and replacing $A(t)$ by $u\to e^{-\lbb t}A(t)(e^{\lbb t}u)$ we can always change $\lbb$ in \eqref{e1.1} to a big enough $\lbb$ which satisfies $\lbb>\nu$.   It should be emphasized that a general existence and uniqueness result for equation \eqref{e1.1} is known only for the special case where $A(t)$ are monotone and demicontinuous operators from $V$ to $V'$, where $(V,V')$ is a pair of reflexive Banach spaces in duality with the Hilbert space $H$ as pivot space, that is,  $V\subset H (\equiv H')\subset V'$ densely and continuously.  If, in addition, for $\a_1\in(0,\9)$, $\a_2,\a_3\in\rr$,
\begin{eqnarray}
{}_{V'}(A(t)u,u)_{V}&\ge&\alpha_1\|u\|^p_V+\alpha_2|u|^2_H,\ \ff u\in V,\label{e1.11}\\[1mm]
\|A(t)u\|_{V'}&\le&\alpha_3\|u\|^{p-1}_V,\ \ff u\in V, \label{e1.12}
\end{eqnarray}where $1<p<\9$, then equation \eqref{e1.1} has under assumptions (i)--(iii)  a unique strong solution $X\in L^p((0,T)\times\ooo;V)$ (see \cite{8}, \cite{13-b}, \cite{9}, \cite{10}). We~noted before that  assumption (i) implies that $A(t,\oo)$ is maximal monotone in $H$ for all $t\in[0,T]$, though not every maximal monotone operator $A(t):D(A(t))\subset H\to H$ has a realization in a convenient pair of spaces $(V,V')$ such that \eqref{e1.8}--\eqref{e1.9} hold.
Though assumptions \eqref{e1.11}--\eqref{e1.12} hold for a large class of stochastic parabolic equations in Sobolev spaces $W^{1,p}(\calo),$ $1\le p<\9$ (see \cite{4}), some other important stochastic PDEs are not covered by this functional scheme. For instance, the variational stochastic differential equations, nonlinear parabolic stochastic equations in $W^{1,1}(\calo)$, in Orlicz-Sobolev spaces on $\calo$  or in $BV(\calo)$ (bounded variation stochastic flows) cannot be treated in this functional setting. As a matter of fact, contrary to what happens for deterministic infinite differential equations, there is no general existence theory for equation \eqref{e1.1} under assumption (i)--(iii). The definition of a convenient concept of a weak solution to be unique and continuous with respect to data is a challenging objective of the existence   theory of the infinite dimensional SDE. In this paper, we introduce such a  solution $X$ for \eqref{e1.1} which is defined as a minimum point of a certain convex functional defined on a suitable space of $H$-valued processes on $(0,T)$. This idea was developed in \cite{6}  for nonlinear operators $A(t):V\to V'$ satisfying condition \eqref{e1.11}--\eqref{e1.12} and is based on the so-called {\it Brezis--Ekeland variational principle} \cite{6}. Such a solution  in the sequel will be called the {\it variational solution} to \eqref{e1.1}. (Along these lines see also \cite{2}, \cite{3}, \cite{4b}, \cite{4}.)

\section{The variational solution to equation \eqref{e1.1}}
\setcounter{equation}{0}

First, we transform  equation \eqref{e1.1} into a random differential equation via the substitution
\begin{equation}\label{e2.1}
X(t)=e^{W(t)}(y(t)+x),\ t\in[0,T],
\end{equation}
which, by It\^o's product rule,

$$dX=e^Wdy+e^W(y+x)dW+\mu e^W(y+x)dt,$$leads to
\begin{equation}\label{e2.2}
\barr{l}
\dd\frac{dy(t)}{dt}+e^{-W(t)}
A(t)(e^{W(t)}(y(t)+x))+(\mu+\lbb)(y(t)+x)\ni 0,\\\hfill t\in(0,T),\vspp
y(0)=0.\earr
\end{equation}
(In the following, we shall omit $\oo$ from the notation $A(t,\oo)$.)

 As a matter of fact, the equivalence between \eqref{e1.1} and \eqref{e2.2} is true only for  a smooth solution $y$ to \eqref{e2.2}, that is, for pathwise absolutely continuous strong solutions to \eqref{e2.2} (see \cite{5}, \cite{5a}). In the sequel, we shall define a generalized (variational) solution for the random Cauchy problem \eqref{e2.2} and will call the corresponding process $X$ defined by \eqref{e2.1} {\it the variational solution to}~\eqref{e1.1}.

We shall treat equation \eqref{e2.2} by the operator method developed in \cite{5a}. Namely,   consider the space $\calh$ of all $H$-valued processes $y:[0,T]\to H$ such that
$$|y|_{\calh}=\(\E\int^s_0|e^{W(t)}y(t)
|^2_{H}dt\)^{\frac12}<\9,$$
which have an $(\calf_t)_{t\ge0}$-adapted version. Here $\E$ denotes the expectation with respect to $\mathbb{P}.$   The space $\calh$ is a Hilbert space with the scalar product
$$\<y,z\>=\E\int^s_0(e^{W(t)}y(t),e^{W(t)}z(t))dt,\ y,z\in\calh.$$

We set $\delta=\frac12\,(\lbb-\nu)$.
Now, consider the operators $\cala:D(\cala)\subset\calh\to\calh$ and $\calb:D(\calb)\subset\calh\to\calh$ defined by
\begin{equation}\label{e2.3}
\barr{rcl}
(\cala y)(t)&\!\!\!=\!\!\!&e^{-W(t)}A(t)(e^{W(t)}(y(t)+x))+\delta (y+x),\ \ff y\in D(\cala ),\\
&&\hfill t\in[0,T],\vspp
D(\cala )&\!\!\!=\!\!\!&\{y\in\calh ;\ e^{W(t)}(y(t)+x)\in D(A(t)),\ \ff t\in[0,T]\mbox{ and }\\
&&\hfill e^{-W}A(e^{W}(y+x))\in\calh \},\earr
\end{equation}
\begin{equation}\label{e2.4}
\barr{rcl}
(\calb y)(t)&\!\!\!=\!\!\!&\dd\frac{dy}{dt}\,(t)+(\mu+\nu+\delta)(y+x),\ \mbox{a.e.}\ t\in(0,T),\ y\in D(\calb ),\vspp
D(\calb )&\!\!\!=\!\!\!&\left\{y\in\calh ;\ y\in W^{1,2}_0([0,T];H),\ \pas,\ \dd\frac{dy}{dt}\in\calh \right\}\!.  \earr
\end{equation}
Here, $W^{1,2}_0([0,T];H)$ denotes the space $\{y\in W^{1,2}([0,T];H);\,y(0)=0\}$, where $W^{1,2}([0,T];H)$ is the Sobolev space $\left\{y\in L^2(0,T;H),\,\frac{dy}{dt}\in L^2(0,T;H)\right\}$. We recall that $W^{1,2}([0,T];H)\subset AC([0,T];H)$,   the space of all $H$-valued absolutely continuous functions on $[0,T]$.

 Then we may rewrite equation \eqref{e2.2} as
\begin{equation}\label{e2.5}
\calb y+\cala y\ni0 .\end{equation}(If $A(t)$ is multivalued,   we replace $A(t)(e^W(y+x))$ in \eqref{e2.3} by $\{\eta(t);\ \eta(t)\in A(t)(e^{W(t)}(y(t)+x)),\mbox{ a.e. }(t,\oo)\in(0,T)\times\ooo\}$.)

Consider the functions $\Phi :\calh \to\ov\rr$ defined by
\begin{equation}\label{e2.6}
\dd\Phi (y)=\E\int^T_0(\vf(t,e^{W(t)}(y(t)+x))+\frac\delta2\,| e^{W(t)}(y(t)+x)|^2_H)dt,\ \ff y\in\calh .\end{equation}It is easily seen that $\Phi $ is convex, lower-semicontinuous and
\begin{equation}\label{e2.7}
\pp\Phi =\cala  .\end{equation}As regards the operator $\calb $,  we have

\begin{lemma}\label{l2.2} For each   $y \in D(\calb )$ we have
\begin{equation}\label{e2.8}
\barr{lcl}
\<\calb y,y\> &=&\dd\frac12\,\E|e^{W(T)}y(T)|^2_H+(\nu+\de)
|y|^2_{\calh }\vspp
&&-\dd\frac12\,\E\int^T_0\sum^\9_{j=1}|e^Wy e_j|^2_H\mu^2_jdt\vspp
&\ge&\dd\frac12\,\E|e^{W(T)}y(T)|^2_H
+\dd\frac\lbb2\,|y|^2_{\calh } .\earr\end{equation}\end{lemma}

\n{\bf Proof.} We have
\begin{equation}\label{e2.9}
\barr{lcl}
\<\calb y, y \> &=&\dd\E\int^T_0\(e^{W(t)}\frac{dy}{dt}\,(t),e^{W(t)
}y(t)\)dt\vspp
&&\dd+\,\E\int^T_0((\mu+\nu+\delta)e^W y,e^W y)dt.\earr\end{equation}Taking into account that
$$d(e^Wy)=e^W \,dy+e^W y\,dW+\mu e^Wy\,dt,\ \ff y\in D(\calb ),$$we get via It\^o's formula that (see \cite{4})
$$\barr{lcl}
\dd\frac12\,d|e^Wy|^2_H&=&\(e^W\,\dd\frac{dy}{dt},e^Wy\)dt+(e^Wy,e^Wy\,dW)+(\mu e^Wy,e^Wy)dt\vspp
&&+\dd\frac12\sum^\9_{j=1}\mu^2_j|e^Wye_j|^2_Hdt.\earr$$Hence
$$\barr{lcl}
\E\dd\int^T_0\(e^W\frac{dy}{dt},e^W y\)dt&=&\dd\frac12\,\E|e^{W(T)}y(T)|^2_H
-\E\dd\int^T_0(\mu e^W y,e^W y)dt\vspp
&&-\dd\frac12\,\E\int^T_0\sum^\9_{j=1}|e^W y e_j|^2_H\mu^2_j\,dt,\earr$$
and so, because $\lbb>\nu$, by \eqref{e1.9}, \eqref{e2.9}, we get \eqref{e2.8}, as claimed.\hf\bk

Consider now the conjugate $\Phi^* :\calh \to\ov\rr$ of functions
$\Phi $, that is,
$$\Phi^* (z)=\sup\{\<z,y\>_{\calh }-\Phi (y);\ y\in\calh \}.$$By \eqref{e2.6}, we see that (see \cite{16-b})
\begin{equation}\label{e2.10}
\Phi^* (z)=\E\int^T_0(\psi^* (t,e^{W(t)}z(t))-(e^{W(t)}z(t),e^{W(t)}x))dt,\end{equation}where $\psi^*$ is the conjugate of the function
\begin{equation}\label{e2.11}
\psi(t,y)=\vf(t,y)+\frac\delta2\,|y|^2_H,\end{equation}
that is,
\begin{equation}\label{e2.11a}
\psi^*(t,v)=\sup\{(v,y)-\vf(t,y)-
\frac\delta2\,|y|^2_H;\ y\in H\}.\end{equation}
We recall  (see \eqref{e1.8}) that
\begin{equation}\label{e2.12}
\Phi (y)+\Phi^* (u)\ge\<y,u\> ,\ \ff y,u\in\calh ,  \end{equation}with equality if and only if $u\in\pp\Phi (y)$. We infer that $y^*$ is  {\it a solution to equation \eqref{e2.5} if and only if}
\begin{equation}\barr{lcl}
y^*&=&{\rm arg}\min\limits_{\hspace*{-6mm}(y,u)\in\cald(\calb)\times\calh}\{\Phi (y)+\Phi^* (u)-\<y,u\> ;\ \calb  y+u=0\}\vspp
&=&{\rm arg}\min\limits_{\hspace*{-6mm}(y,u)\in\cald(\calb)\times\calh}\{\Phi (y)+\Phi^* (u)+\<\calb y,y\> ;\ \calb y+u=0\}\earr\label{e2.13}
\end{equation}and
\begin{equation}\label{e2.14}
\Phi (y^*)+\Phi^* (u^*)+\<\calb y^*,y^*\> =0.\end{equation}
Taking into account \eqref{e2.10} and  recalling \eqref{e2.6}, \eqref{e2.8}, we have
 \begin{equation}\label{e2.15}
\hspace*{-4mm}\barr{ll}
y^*={\rm arg}\min\limits_{\hspace*{-6mm}(y,u)\in\cald(\calb)\times\calh}\Bigg\{\!\!\!&\E\dd\int^T_0 \Big(\vf(t,e^{W(t)}(y(t)+x))
+\dd\frac\delta2\,|e^{W(t)}(y(t)+x)|^2_H\vsp &+\,\psi^*(t,e^{W(t)}u(t))-
(e^{W(t)}u(t),e^{W(t)}x)\\
&+\eta(e^{W(t)}y(t))\Big)dt
+\dd\frac12\,\E|e^{W(T)}y(T)|^2_H;\ \calb y+u=0\Bigg\},\earr\hspace*{-15mm}
\end{equation}
where
\begin{equation}\label{e2.16}
\eta(z)=(\nu+\delta)|z|^2_H-\frac12\sum^\9_{j=1}|z e_j|^2_H\mu^2_j.\end{equation}
We note also that, by It\^o's product rule, we have, for $u\in\calh,$ $y\in\cald(\calb)$,
$$
\barr{l}
-\E\dd\int^T_0(e^{W(t)}u(t),e^{W(t)}x)dt\\
\qquad =\E\dd\int^T_0\bigg(e^Wx,e^W\(\frac{dy}{dt}
+(\mu+\nu+\de)(y+x)\)\bigg)dt\\
\qquad
=\E\dd\int^T_0(e^Wx,d(e^Wy))
-\int^T_0(e^Wx,\mu e^Wy
-(\mu+\nu+\de)(y+x)e^W)dt\\
\qquad=\E\dd\int^T_0(e^Wx,(\mu+\delta+\nu)(y+x)e^W)dt\\
\qquad
+\E\dd\int^T_0d(e^Wx,e^Wy)-\int^T_0
\((e^Wy,e^W(1+\mu) x)-(\mu e^Wy,e^Wx)\)dt\\
\qquad
=\E(e^{W(T)}x,e^{W(T)}y(T))
-\dd\int^T_0\((e^Wy, e^W(1+\mu) x)
-(\mu e^Wy,e^Wx)\)dt\\
\qquad
+\E\dd\int^T_0(e^Wx,(\mu+\nu+\de)
(y+x)e^W)dt\earr$$
$$\barr{l}
=\E(e^{W(T)}x,e^{W(T)}y(T))
+\E\dd\int^T_0(e^Wx,
((\mu+\nu+\de)(y+x)-\mu y)
e^W)dt\\
\qquad-\dd\int^T_0(e^Wy,e^W(1+\mu)x)dt\\
\qquad
=\E\dd\int^T_0(e^{W}
((\nu+\delta)(y +x)+\mu x),
e^{W }x)dt
+\E(e^{W(T)}y(T),e^{W(T)}x)\\
\qquad-\dd\int^T_0(e^Wy,e^W(1+\mu)x)dt.\earr$$

Let $\calh_0$ denote the set of all $u\in L^2([0,T]\times\ooo;H)$ which have an $(\calf_t)_{t\ge0}$-adapted version.
We set, for $y\in\calh,\ u\in\calh_0,$
\begin{eqnarray}
\!\!\!G_1(y)\!\!\!&{=}
& \E\int^T_0 \vf(t,e^{W(t)}(y(t)+x))dt\label{e2.17}\\ &&+\,\E\dd\int^T_0((e^{W(t)}
((\nu+\delta)(y(t)+x)
+\mu x),e^{W(t)}x)\nonumber\\
&&+\,\dd\frac\delta2\,|e^{W(t)}(y(t)+x)|^2_H
+\eta(e^{W(t)}y(t)))dt\nonumber\\[1mm]
&&+\dd\frac12\,\E|e^{W(T)}y(T)|^2_H
+\E(e^{W(T)}y(T),e^{W(T)}x)\nonumber\\
&&-\E\dd\int^T_0(e^{W(t)}y(t),e^{W(t)}
(1+\mu)x)dt,\nonumber\\[1mm]
\!\!\!G_2(u)\!\!\!&{=}&\!\!\!\E\int^T_0 \psi^*(t,u(t))  dt,\quad\label{e2.18}
\end{eqnarray}where $\psi^*$ is given by \eqref{e2.11a}.

By \eqref{e2.15} it follows that $y^*$ is a solution to equation \eqref{e2.5} {\it if~and only if}
\begin{equation}\label{e2.19}
y^*={\rm arg}\min\limits_{(y,u)\in\cald(\calb)\times\calh_0}\{G_1(y)+G_2(u);\ e^W \calb y+u=0\}\end{equation}{\it and}
\begin{equation}\label{e2.20}
G_1(y^*)+G_2(u^*)=0.\end{equation}
It should be said, however, that under our assumptions the convex minimization problem \eqref{e2.19} might have no solution $(y^*,u^*)$ because, in general, $G_2$ is not coercive on the space $\calh$. ($G_2$ is, however, coercive if $\vf$ is bounded on bounded sets of $H$. But such a condition is too restrictive for applications to PDEs.) So, we are led to replace \eqref{e2.19} by a relaxed optimization problem to be defined below.

Let
\begin{equation}\label{e2.20a}
\calx = L^2(\ooo;(W^{1,2}([0,T];H))'), \end{equation}
where $(W^{1,2}([0,T];H))'$ is the dual space of $W^{1,2}([0,T];H)$.

Define the operator $\wt \calb : \calh \times L^2(\ooo;H) \to \calx$ by

\begin{equation}\label{e2.21}
\barr{ll}
(\wt\calb(y,y_1))(\theta) \!\!\!&=
 \E(e^{W(T)}y_1,\theta(T))
 +\E\dd\int^T_0((\nu+\delta)
 (y(t)+x)\vsp
 &
 +\mu x)e^{W(t)},\theta(t))dt
 -\E\dd\int^T_0\(e^{W(t)}y(t),
 \frac{d\theta}{dt}\,(t)\)dt,\ \vsp
 &\hfill\ff\theta\in L^2(\ooo;W^{1,2}([0,T];H)).\earr
 \end{equation}
We note that $y_1(\oo) \in H$ can be viewed as the trace   of $y(\oo)$ at $t=T$.

Indeed, if $y\in \cald(\calb)$, we have via It\^o's formula
$$\barr{lcl}
\E\dd\int^T_0(e^W\calb y,\theta)dt\!\!\!
&=\!\!\!&\E\(\dd\int^T_0(d(e^Wy),\theta)
-\int^T_0(e^W\mu y, \theta)dt\)\vsp
&&
+\E\dd\int^T_0(e^W(\mu+\nu+\delta)
(y+x),\theta)dt\vsp
&=\!\!\!&\E(e^{W(T)}y(T),\theta(T))\vsp
&&
+\E\dd\int^T_0(e^W((y+x)(\nu+\delta)+\mu x),\theta)dt\vsp
&&-\E\dd\int^T_0\(e^Wy,
\frac{d\theta}{dt}\)dt,\
\ff\theta\in L^2(\ooo;W^{1,2}([0,T];H)).\earr$$
This means that $\wt\calb(y,y(T))=e^W\calb y,$ $\ff y\in \cald(\calb)$. We set
\begin{equation}\label{e2.23-z}
\barr{ll}
\wt G_1(y,y_1)=\!\!\!
&\E\dd\int^T_0\vf
(t,e^{W(t)}(y(t)+x))dt\vsp
&
+\E\dd\int^T_0\bigg((e^{W(t)}
((\nu+\delta)(y(t)+x)+\mu x),e^{W(t)}x) \vsp
&
+\dd\frac\delta2\,
|e^{W(t)}(y(t)+x)|^2_H
+\eta(e^{W(t)}y(t))\bigg)dt\vsp
&-\E\dd\int^T_0
(e^{W(t)}y(t),e^{W(t)}
(1+\mu)x)dt
+\dd\frac12\,\E|e^{W(T)}y_1|^2_H\vsp
&+\E(e^{W(T)}y_1,e^{W(T)}x),\
 \  \ff (y,y_1)\in\calh\times L^2(\ooo;H)\earr\end{equation}
and note that $\wt G_1(y;y(T))=G_1(y),\ \ff y\in\cald(\calb)$.

We note also that, if  $y_n\in\cald(\calb)$ such that $y_n\to y$ weakly in $\calh$ and $y_n(T)\to y_1$ weakly in $L^2(\ooo;H)$, then
\begin{equation}\label{e2.23-zz}
e^W\calb y_n\to\wt\calb(y,y_1)\mbox{ weakly in }\calx.\end{equation}

Let $\ov G:\calh\times L^2(\ooo;H)\times \calx\to\ov\rr$ be the lower semicontinuous closure of the function $G(y,y_1,u)=\wt G_1(y,y_1)+G_2(u)$ in $\calh\times L^2(\ooo;H)\times\calx$, on the set $\{(y,y_1,u)\in\calh\times L^2(\ooo;H)\times\calx;\ e^W\calb y+u=0\},$ that is,
\begin{equation}\label{e2.22}
\barr{r}
\ov G(y,y_1,u)=\liminf
\{G(z,z(T),u);\,z(T)\to y_1\mbox{ in }L^2(\ooo;H),\vsp
z\in\cald(\calb),\,(z,v)\to(y,u)\mbox{ in }\calh\times\calx;\,
e^W \calb z+v=0\}. \earr \end{equation}
(Here and everywhere in the following, by $\to$ we mean weak
convergence.)

Taking into account that the function $\wt G_1$ is convex and lower semicontinuous in $\calh\times L^2(\ooo;H)$, we have by \eqref{e2.22}
\begin{equation}\label{e2.24-b}
\barr{r}  \ov G(y,y_1,u){=}\wt G_1(y,y_1){+}\liminf
\{G_2(v); (z,v){\to}(y,u)\mbox{ in }\calh{\times}\calx,\vsp  e^W\calb z{+}v{=}0\}.\earr
\end{equation}

Now, we relax \eqref{e2.19} to the convex minimization problem
$${\rm Min}\{\ov G(y,y_1,u);\, \widetilde\calb(y,y_1)+u=0;\ (y,y_1,u)\in\calh\times L^2(\ooo;H)\times\calx\}.\leqno{\rm(P)}$$
We have
\begin{theorem}\label{t2.1} Let $x\in H$. Then problem {\rm(P)} has a unique solution $(y^*,y^*_1,u^*)$ $\in\!\calh{\times } L^2(\ooo;H){\times}\calx,$ with $u^*{=}-\wt\calb(y^*,y^*_1)$.
Moreover, $\vf(\cdot,e^{W}(y^*+x))\!\in\! L^1((0,T){\times}\ooo).$\end{theorem}
\n{\bf Proof.} Let $m$ be the infimum in (P)  and let $(y_n,u_n)\in \cald(\calb)\times\calh$ be such that
\begin{eqnarray}
&m\le G(y_n,y_n(T),u_n)\le m+\dd\frac1n,\ \ \ff n\in\mathbb{N},\label{e2.23}\\
&e^W\calb y_n+u_n=0.\label{e2.24}
\end{eqnarray}
Since, by assumption (iii), for some $C_1,C_2\in]0,\9[$,
  $$\wt G_1(y_n,y_n(T))\ge C_1(|y_n|^2_\calh+\E|e^{W(T)}y_n(T)|^2_H)-C_2,$$
    we have along a subsequence
$$y_n\longrightarrow y^*\mbox{ weakly in }\calh,\ \ y_n(T)\to y^*_1\mbox{ weakly in }L^2(\ooo;H),$$and so, by \eqref{e2.23-zz}, we have
$$u_n\longrightarrow u^*=-\wt\calb (y^*,y^*_1)\mbox{ weakly in }\calx.$$

As $\ov G$ is  weakly  lower semicontinuous on $\calh\times L^2(\ooo;H)\times\calx$,   we see by \eqref{e2.23} that
$$\ov G(y^*,y^*_1,u^*)=m,$$as claimed.
The  uniqueness of $(y^*,y^*_1,u^*)$ is immediate because the function $\ov G(\cdot,\cdot,u)$ is  strictly convex on~$\calh\times L^2(\ooo;H)$ for all $u\in\calx$.\hf

\begin{definition}\label{d2.2}\rm A pair $(y^*,y^*_1)$ such that  $(y^*,y^*_1,u^*)\in\calh\times L^2(\ooo;H)\times\calx$,  $u^*=-\wt B(y^*,y^*_1)$, is a solution to problem (P), is called {\it the variational solution to equation} \eqref{e2.2}, and $X^*=e^W(y^*+x)$ is called the {\it variational solution to equation} \eqref{e1.1}.\end{definition}

The variational solution $X^*:(0,T)\to H$ is an $(\calf_t)_{t\ge0}$-adapted process.

 Theorem \ref{t2.1} can be rephrased as:
 \begin{theorem}\label{t2.4-b} Under hypotheses {\rm(i)--(iii)}, equation  \eqref{e1.1} has a unique va\-ria\-tio\-nal solution $X^*\in L^2((0,T)\times\ooo;H)$ with $\vf(t,X^*)\in L^1((0,T)\times\ooo).$
 \end{theorem}
 It should be noted that $y^*$ and $X^*$, as well, are not pathwise continuous on $[0,T]$. As seen later on, this happens, however, in some specific cases with respect to a weaker topology.

In the next section, we shall see how problem (P) looks like in a few important examples of stochastic PDEs.
\begin{remark}\rm The above formulation of the variational solution $X^*$ is strongly dependent on the subdifferential form \eqref{e1.3} of the operator $A(t)$.  The extension of the above technique to a general maximal monotone function $A(t):H\to H$ remains to be done using the Fitzpatrick formalism (see~\cite{11}).\end{remark}
\section{Nonlinear parabolic stochastic differential equations}
\setcounter{equation}{0}
We consider here the stochastic differential equation
\begin{equation}\label{e3.1}
\barr{l}
dX-{\rm div}_\xi(a(t,\nabla X))dt+\lbb X\,dt=X\,dW\mbox{ in }(0,T)\times\calo,\vspp
X=0\mbox{ on }(0,T)\times\pp\calo,\vspp
X(0,\xi)=x(\xi),\ \xi\in\calo\subset\rr^d,\earr\end{equation}
where $x\in H$,  $W$ is the Wiener process \eqref{e1.2} in $H=L^2(\calo)$, $\calo$ is a bounded and open subset of $\rr^d$ with smooth boundary $\pp\calo$, and  $a:(0,T)\times\rr^d\to\rr^d$ is a nonlinear mapping of the form
\begin{equation}\label{e3.2}
a(t,z)=\pp_z j(t,z),\ \ff z\in\rr^d,\ t\in[0,T],\end{equation}
where $j:(0,T)\times\rr^d\to\rr$ is measurable, convex, lower semicontinuous in $z$ and
\begin{eqnarray}
\lim_{|z|\to\9}\frac{j(t,z)}{|z|}&=&+\9,\ \ t\in[0,T],\label{e3.3}\\[1mm]
\lim_{|v|\to\9}\frac{j^*(t,v)}{|v|}&=&+\9,\ \ t\in[0,T],\label{e3.4}
\end{eqnarray}uniformly with respect to $t\in[0,T].$

We note that, if the function  $(t,y)\to j(t,y)$ is bounded on bounded subsets of $[0,T]\times\rr^d$, then \eqref{e3.4} automatically holds by the conjugacy formula \eqref{e1.8}, that is,
$$j^*(t,v)\ge v\cdot z-j(t,z),\ \ff v,z\in \rr^d,\ t\in[0,T].$$

It should be noted that equation \eqref{e3.1} cannot be treated in the functional setting \eqref{e1.11}-\eqref{e1.12} which require polynomial growth and boundedness  for $j(t,\cdot)$, while assumptions \eqref{e3.3}--\eqref{e3.4} allow nonlinear diffusions $a$ with slow growth  to $+\9$ as well as superlinear growth  of the form
$$a(t,z)=a_0\exp(a_1|z|^p{\rm sgn}\ z).$$

We   note also  that assumptions \eqref{e3.2}--\eqref{e3.4} do not preclude multi\-valued mappings $a$.  Such an example is
$$\barr{rcl}
j(t,z)&\equiv&|z|(\log(|z|+1)),\vspp
a(t,z)&=&\(\log(|z|+1)+\dd\frac1{|z|+1}\){\rm sign}\ z,\ \ \ff z\in\rr^d.\earr$$

By \eqref{e2.1}, one reduces equation \eqref{e3.1} to the random parabolic differential equation
\begin{equation}\label{e3.5}
\barr{l}
\dd\frac{\pp y}{\pp t}-e^{-W}{\rm div}_\xi a(t,\nabla(e^W(y+x))+(\lbb+\mu)(y+x))=0\\
\hfill\mbox{ in }(0,T)\times\calo,\vspp
y=0\mbox{ on }(0,T)\times\pp\calo,\vspp
y(0,\xi)=0,\ \ \xi\in\calo.\earr\end{equation}We are under the conditions of Section 2, where
$$\barr{rcl}
 H&=&L^2(\calo), \vsp
 A(t)y&=&-{\rm div}_\xi a(t,\nabla y),\vsp
  D(A(t))&=&\{y\in W^{1,1}_0(\calo);\ {\rm div}_\xi a(t,\nabla y)\in L^2(\calo)\}\vsp
  \vf(t,y)&=&\dd\int_\calo j(t,\nabla y(\xi))d\xi.\earr$$

By \eqref{e2.11a}, we have
\begin{equation}\label{e3.6}
\psi^*(t,v)=\int_\calo(a(t,\nabla z)\cdot\nabla z-j(t,\nabla z)+\frac\delta2\, z^2)d\xi,\ \ff v\in L^2(\calo),\end{equation}where $z$ is the solution to the equation
\begin{equation}\label{e3.7}
\barr{rcll}
-{\rm div}\ a(t,\nabla z)+ \delta z&=&v&\mbox{ in }\calo,\vspp
z&=&0&\mbox{ on }\pp\calo,\earr
\end{equation}
or, equivalently,
\begin{equation}\label{e3-7prim}
z={\rm arg}\min_{\tilde z\in W^{1,1}_0(\calo)}
\left\{
\int_\calo j(t,\nabla \tilde z)d\xi
-\int_\calo v\tilde z\,d\xi+\frac\delta 2
\int_\calo \tilde z^2d\xi\right\}.
\end{equation}
 By \eqref{e3.3}, it follows that \eqref{e3-7prim} has, for each $v\in L^2(\calo)$ and $t\in[0, T]$, a~unique solution $z\in W^{1,1}_0(\calo)$. In fact, as easily seen, by condition \eqref{e3.3} it follows that the functional arising in the right side part of \eqref{e3-7prim} is convex, lower semicontinuous and coercive on $W^{1,1}_0(\calo)$. By  \eqref{e2.23-z}), we have
\begin{equation}\label{e3.8}
 \barr{ll}
\wt G_1(y,y_1)=\!\!\!&\dd
\E \dd\int^T_0   \int_\calo (j(t,\nabla(e^{W(t)}(y(t)+x)))
+\dd\frac\delta2\,|e^{W(t)}(y(t)+x)|^2_H\vspp
&+\,\dd e^{W(t)}
((\nu+\delta)(y(t)+x)+\mu x)e^{W(t)}x)\,d\xi\,dt\vsp
&-\E\dd\int^T_0(e^{W(t)}
y(t),e^{W(t)}(1+\mu)x)dt\vsp
&+ \dd\E\int^T_0 \eta(e^{W(t)}y(t))dt
+ \dd\frac12\,\E \int_\calo |e^{W(T)}y_1(\xi)|^2 d\xi\vsp
&+\E(e^{W(T)}y_1,e^{W(T)}x),
\,(y,y_1)\in\calh\times L^2(\ooo;H),\earr\hspace*{-6mm}
\end{equation}where (see \eqref{e2.16})
\begin{equation}\label{e3.9}
\eta(z)=(\nu+\delta)\int_\calo|z|^2d\xi-\frac12\sum^\9_{j=1}\mu^2_j\int_\calo|z e_j|^2d\xi.\end{equation}
By \eqref{e2.18} and \eqref{e3.6}--\eqref{e3.7}, we also have \begin{equation}\label{e3.10}
\barr{ll}
G_2(u)=\!\!\!
&\dd\E\int^T_0   \int_\calo(a(t,\nabla z(t,\xi))\cdot\nabla z(t,\xi))\\
&-\dd j(t,\nabla z(t,\xi))+\dd\frac\delta 2\, z^2(t,\xi)d\xi\, dt,\ u\in\calh,\earr\end{equation}
where $z(t,\oo)\in  W^{1,1}_0(\calo )$ for $dt\otimes\mathbb{P}$-a.e.,   $(t,\oo)\in (0,T)\times\ooo $, is given by (see~\eqref{e3.7})
\begin{equation}\label{e3.11}
\barr{rcll}-{\rm div}\ a(t,\nabla z)+ \delta z&=&u&\mbox{ in }\calo,\vspp
z&=&0&\mbox{ on }\pp\calo.\earr\end{equation}
Taking into account that $a(t,\nabla z)\cdot\nabla z\ge j(t,\nabla z)-j(t,0)$, we see by \eqref{e3.3} and \eqref{e3.11} that
$$z\in L^1((0,T)\times\ooo;W^{1,1}(\calo))\cap L^2((0,T)\times\calo\times\ooo).$$
Recalling \eqref{e1.7}--\eqref{e1.8}, we have
$$a(t,\nabla z)\cdot\nabla z-j(t,\nabla z)=j^*(t,a(t,\nabla z))\mbox{ a.e. in }(0,T)\times\calo,$$and this yields
\begin{equation}\label{e3.12}
\barr{ll}
G_2(u)= \dd\E\int^T_0
\int_\calo (j^*(t,a(t,\nabla z(t,\xi)))
 +\dd\frac\delta2\, z^2(t,\xi) )d\xi\,dt.\earr\end{equation}
By \eqref{e3.4}, it follows via the Dunford--Pettis weak compactness theorem in $L^1$ that every level set
$$\left\{v;\E\int^T_0\int_\calo j^*(t,v(t,\xi))dx\,d\xi\le M\right\},\ M>0,$$is weakly compact in the space $
L^1((0,T)\times\calo\times\ooo). $ By \eqref{e3.11} and \eqref{e3.12}, we see that, if
 $G_2(u_n)\le M,$ where $ \{u_n\}\subset L^2((0,T)\times\calo\times\ooo\}$ and $z_n$ is the solution to \eqref{e3.11} with $u_n$ replacing $u$, then, by the Dunford--Pettis  theorem, the sequence $\{a(t,\nabla z_n)\}$ is weakly compact in $L^1 ((0,T)\times\calo\times\ooo\}$. Hence $\{u_n\}$ is weakly compact in $L^1 ((0,T)\times\ooo;W^{-1,\9}(\calo))$.

 By \eqref{e3.12}, it follows also that $\{z_n\}$ is weakly compact in $L^2((0,T)\times\calo\times\ooo).$\break
    By \eqref{e2.22}, this means that, if $x\!\in\! L^2(\calo)$, then, for $(y,y_1,u)\!\in \!\calh{\times} L^2(\ooo;H){\times}\calx$,
\begin{equation}\label{e3.14}
\barr{ll}
\ov G (y,y_1,u)\vsp
= \wt G_1(y,y_1)+
 \dd\E\int^T_0\!\!\! \int_\calo\!\bigg(j^*(t,a(t,\nabla z(t,\xi)) +\frac\delta2\, z^2(t,\xi))\!\bigg)d\xi\,dt, \earr \end{equation}
where   $z\in L^1((0,T)\times\ooo;W^{1,-1}_0(\calo))\cap L^2((0,T)\times\calo\times\ooo)$ is the solution to \eqref{e3.11}.

Let $(y_n,u_n)\in\calh\times\calh$ be such that $e^W\calb y_n+u_n=0$ and $(y_n,u_n)\to(y,u)$ in $\calh\times\calx,$ $y_n(T)\to y_1)$ in $L^2(\ooo;H)$.
Since $\sup\limits_n\{G_1(y_n)\}<\9$, by \eqref{e3.3} and \eqref{e3.8}, it follows also that $\{\nabla(e^W(y_n+x)\}$ is weakly compact in $L^1 ((0,T)\times\calo\times\ooo)$, and so   $e^W(y+x)\in L^1 ((0,T)\times\calo;W^{1,1}_0(\calo)).$ Moreover,   it follows that $\left\{\frac{dy_n}{dt}\right\}$ is weakly compact in $L^1((0,T)\times\ooo;W^{-1,\9}(\calo))$, and so $\frac{dy}{dt}\in L^1((0,T)\times\ooo;W^{-1,\9}(\calo)).
$ This implies that the equation  $\wt\calb(y^*,y^*_1)+u^*=0$ reduces to
$$\barr{l}
e^W\,\dd\frac{dy^*}{dt}+e^W
(\mu+\nu+\delta)(y^*+x)
+u^*=0\mbox{ in }\cald'(0,T),\ \pas,\vsp
y^*(0)=0,\quad y^*(T)=y^*_1.\earr$$
Hence, if $D(\ov G_1)=\{(y,y_1,u);\,\ov G_1(y,y_1,u)<\9\}$, then we have
$$\barr{l}
D(\ov G_1)\subset\{(y,y_1,u)\in\calh\times L^2(\ooo;H)\times\calh;\,e^Wy\in L^1((0,T)\times\ooo; W^{1,1}_0(\calo));\vsp
\dd\frac{dy}{dt}\in L^1((0,T)\times\ooo;W^{-1,\9}(\calo));\,u\in L^1((0,T)\times\ooo;W^{-1,\9}(\calo)),
\,y_1=y(T)\}.\earr$$
This means that, in this case, problem (P) can be rewritten as
\begin{equation}\label{e3.15}
\barr{ll}
{\rm Min}\Big\{\!\!\!&\ov G(y,y(T),u);\,y\in L^2((0,T)\times\calo\times\ooo)\cap \calh,\\
& e^W(y+x)\in L^1((0,T)\times\ooo;W^{1,1}_0(\calo)),\vsp
&\dd\frac{dy}{dt}\in L^1((0,T)\times\ooo;W^{-1,\9}(\calo)),\vsp & u\in L^1((0,T)\times\ooo;W^{-1\9}(\calo))\cap\calx;\vspp
&\mbox{subject to }\vsp
&\dd\frac{dy}{dt}+(\mu+\nu+\delta)(y+x)+e^{-W} u=0\mbox{ on }(0,T);\,y(0)=0\Big\},\earr
\hspace*{-10mm}\end{equation}
where $\ov G_1$ is defined by \eqref{e3.14}.
By Theorem \ref{t2.1}, there is a unique solution $(y^*,u^*)$ to \eqref{e3.15}. Taking into account that $ u^*\in L^1((0,T)\times\ooo; W^{-1,\9}(\calo))$ and that
$$y^*(t)=-\int^t_0 e^{-W}u^*(s)ds-\int^t_0(\mu+\nu+\delta)(y^*(s)+x)ds,\ \ff t\in(0,T),$$ we infer that the process $t\to y^*(t)$ in pathwise $W^{-1,\9}(\calo)$ continuous on$(0,T)$. By Theorem \ref{t2.4-b}, we have, therefore,

\begin{theorem}\label{t3.1} Assume that $x\in L^2(\calo)$ and that conditions \eqref{e3.2}--\eqref{e3.4} hold. Then, equation \eqref{e3.1} has a unique variational solution
\begin{equation}\label{e3.16}
X^*\in L^2((0,T)\times\calo\times\ooo),\ e^WX^*\in L^1((0,T)\times\ooo;W^{1,1}_0(\calo)).\end{equation}
Moreover, the process $t\to X^*(t)$ is $(\calf_t)_{t\ge0}$-adapted and pathwise $W^{-1,\9}(\calo)$-valued  continuous on $(0,T)$.\end{theorem}

\paragraph{The total variation flow}\ \bk

\n The stochastic differential equation
\begin{equation}\label{e3.17}
\barr{l}
dX-{\rm div}\(\dd\frac{\nabla X}{|\nabla X|_d}\)dt+\lbb X\,dt=X\,dW\mbox{ in }(0,T)\times\calo,\vspp
X(0)=x\mbox{ in }\calo,\vspp
X=0\mbox{ on }(0,T)\times\pp\calo\earr
\end{equation}
with $x\in L^2(\calo)$ is the equation of stochastic variational flow in $\calo\subset\rr^d$, $1\le d\le 3$. The existence and uniqueness of a generalized solution to \eqref{e3.17} $X:[0,T]\to BV(\calo)$ was established in \cite{5} by using some specific approximation techniques. We shall treat now equation \eqref{e3.17} in the framework of variational solution developed above in the space $H=L^2(\calo)$ with the norm $|\cdot|_H=|\cdot|_2$ and the scalar product $(\cdot,\cdot)$, and $\vf:L^2(\calo)\to\ov\rr$ defined by
$$\vf(y)=\left\{\barr{ll}
\|Dy\|+\dd\int_{\pp\calo}|\g_0(y)|d\calh^{d-1},&y\in BV(\calo)\setminus L^2(\calo),\vspp
+\9&\mbox{otherwise}.\earr\right.$$Here, $BV(\calo)$ is the space of functions with bounded variation and $\|Dy\|$ is the total variation of $y\in BV(\calo).$ (See, e.g., \cite{5}.)
Then, with the notations of Section 2, we have  $Ay=\pp\vf(y)$, where $\pp\vf:L^2(\calo)\to L^2(\calo)$ is the subdifferential of $\vf$ and (see \eqref{e2.2}, \eqref{e2.17})
\begin{equation}\label{e3.18}
\barr{l}
\dd\frac{\pp y}{\pp t}+e^{-W}A(e^W(y+x))+\mu(y+x)=0\mbox{ in }(0,T)\times\calo,\vspp
y(0,\xi)=0,\ \ \xi\in\calo,\vspp
y=0\mbox{ on }(0,T)\times\pp\calo.\earr\end{equation}
The function $\wt G_1$ is given, in this case, by

\begin{equation}\label{e3.19}
\barr{ll}
\wt G_1(y,y_1)=\!\!\!&\E\dd\int^T_0 \bigg(\vf(e^{W(t)}(y(t)+x))+
\dd\frac\delta2\,|e^{W(t)}(y(t)+x)|^2_2\vspp
&+\dd (e^{W(t)}((\nu+\delta)(y(t)+x)+
\mu  x),e^{W(t)}x)\bigg)dt  \vspp
&-\E\dd\int^T_0(e^{W(t)}y(t),e^{W(t)}(1+\mu)x)dt\vspp
&+ \E\dd\int^T_0 \eta(e^{W(t)}y(t))dt+
\dd\frac12\,\E|e^{W(T)}y_1|^2_2\vsp &+\E(e^{W(T)}y_1,e^{W(T)}x),\,
(y,y_1)\in\calh\times L^2(\ooo;H),\earr\end{equation}where $\eta$  is given by \eqref{e3.9}.
We have also (see \eqref{e2.11}, \eqref{e2.11a}, \eqref{e3.6})
$$\barr{rcll}
\psi(y)&=&\vf(y)+\dd\frac\delta2\,|y|^2_2,& \ff y\in D(\vf),\vspp
\psi^*(v)&=&(v,\theta)-\vf(\theta)
-\dd\frac\delta2\,|\theta|^2_2,
&v\in\pp\vf(\theta)+\delta \theta.\earr$$
Hence,
$$\barr{lcl}
\psi^*(v)&=&\dd\frac\delta2\,|\theta|^2
+(\pp\vf(z),\theta)-\vf(\theta)\vsp
&=&\dd\frac\delta2\,|(\delta I+\pp\vf)^{-1}v|^2_2+\vf^*(v-(\delta I+\pp\vf)^{-1}v)\earr$$
and, therefore, by \eqref{e2.18},
$$G_2(u)=\E\int^T_0
\(\dd\frac\delta2\,|(\delta I+\pp\vf)^{-1}(u)|^2_2
+\vf^*(u-(\delta  I+\pp\vf)^{-1}( u))\)dt,$$where $\vf^*:L^2(\calo)\to\ov\rr$ is the conjugate of the function $\vf$.    This yields
\begin{equation}\label{e3.20}
\barr{ll}
\ov G(y,y_1,u)=\!\!\!
&G_1(y,y_1)+
\dd\liminf_{^{(z,v)\to (y,u)}_{\  {\rm in\ }  \calh\times\calx}}
\Big\{\E\int^T_0\Big(
\dd\frac\delta2\,
|(I+\pp\vf)^{-1}(v(t))|^2_2\vspp
&\hfill\dd+\vf^*(v-( \delta I+\pp\vf)^{-1}(v(t)))\Big)dt
,\,e^W\calb z+v=0\Big\}, \earr\end{equation}
where the space $\calx$ is defined by \eqref{e2.20a}.

By definition, the solution $(y^*,y^*_1)$ to the minimization problem
\begin{equation}\label{e3.21}
{\rm Min}\{\ov G(y,y_1,u) ;\ \calb (y,y_1)+u=0,\ (y,y_1,u)\in\calh\times L^2(\ooo;H)\times\calx\}\end{equation}  is the {\it variational solution} to the random differential equation \eqref{e3.18}.
\newpage

Denote by $V^*$ the dual of the space $V=BV(\calo)\cap L^2(\calo)$. We note that $\vf^*$ can be extended as a convex lower semicontinuous convex function on $F^*$, and we also have
$$\frac{\vf^*(u)}{\|u\|_{V^*}}
\longrightarrow
+\9\mbox{\ \ as\ \ }\|u\|_{V^*}\longrightarrow +\9.$$
Then, if $(z_n,y_n)\in\calh\times\calh$ is convergent to $(y,u)\in\calh\times\calx$, it follows by the Dunford-Pettis compactness criterium (see \cite{11-a}) that $\{v_n\}$ is weakly compact in $L^1((0,T)\times\ooo;V^*)$. This implies that
$$D(\ov G)\subset L^1((0,T)\times\ooo;BV(\calo))\times L^2(\ooo;H)\times L^1((0,T)\times\ooo;V^*),$$and so, in particular, it follows that
$$y\in W^{1,1}([0,T];V^*),\ \pas$$We have, therefore,
\begin{theorem}\label{t3.2} Let $x\in BV(\calo)\cap L^2(\calo)$. Then equation \eqref{e3.17} has a unique varia\-tio\-nal solution $X=e^W(y+x)$ which is $V^*$-valued pathwise continuous and satisfies
\begin{eqnarray}
 \vf(X)&\in & L^1((0,T)\times\ooo),\label{e3.22}\\[1mm]
 X&\in & L^2((0,T)\times\calo\times\ooo),\ AX\in L^1((0,T)\times\ooo;V^*),\label{e3.23}\\[1mm]
 e^{-W}X&\in&W^{1,1}([0,T];V^*),\ \pas\label{e3.23-az}
\end{eqnarray}
\end{theorem}

In \cite{5}, it was proved the existence and uniqueness of a generalized solution $X$, also called {\it the variational solution}, which was obtained as limit $X^*=\lim\limits_{\vp\to0}X_\vp$ in $L^2(\ooo;C((0,T);L^2(\calo)))$, where $X_\vp$ is the solution to the approximating equation
\begin{equation}\label{e3.24}
\barr{l}
dX_\vp-{\rm div}\ a_\vp(\nabla X_\vp)dt+\lbb X_\vp=X_\vp dW\mbox{ in }(0,T)\times\calo,\vspp
X_\vp(0)=x,\ \ \ X_\vp=0\mbox{ on }(0,T)\times\calo,\earr\end{equation}where $a_\vp=\nabla j_\vp$ and $j_\vp$ is the Moreau--Yosida approximation of the function $r\to|r|_d$. Since, as strong solution to \eqref{e3.24}, $X_\vp$ is also a variational solution to this equation in sense of Definition \ref{d2.2}, it is clear by the structural stability of convex minimization problems that, for $\vp\to0$, we have also $X_\vp\to X$, where $X$ is the variational solution given by Theorem \ref{t3.2}.

We may infer, therefore, that {\it the function $X$ given by Theorem {\rm\ref{t3.2}} is just the generalized solution of \eqref{e3.17} given by Theorem  {\rm3.1}} in \cite{5}. In particular, this implies that $X$ is $L^2(\calo)$-valued pathwise continuous.

In \cite{3-a}, it is de\-ve\-loped a direct variational approach to \eqref{e3.17}, which leads via first order conditions of optimality  to sharper results. (On these lines, see also \cite{10-a}.)

\paragraph{Stochastic porous media equations}\ \bk

\n Consider the equation
\begin{equation}\label{e3.25}
\barr{l}
dX-\Delta\beta(X)dt+\lbb X\,dt=X\,dW\mbox{ in }(0,T)\times\calo,\vspp
X=0\mbox{ on }(0,T)\times\pp\calo,\vspp
X(0,\xi)=x(\xi),\ \ \xi\in\calo,\earr
\end{equation}
where $\calo$ is a bounded and open domain of $\rr^d$, $d\ge1$, $\lbb>0$, $W$ is a Wiener process in $H=H^{-1}(\calo)$ of the form \eqref{e1.2} and $\b$ is a continuous and monotonically nondecreasing function such that $\beta(0)=0$ and
\begin{equation}\label{e3.26}
\lim_{|r|\to\9}\ \frac{j(r)}{|r|}=+\9 .
\end{equation}
In this case, 
$$\barr{rcl}
H&=&H^{-1}(\calo) ,\vsp Ay&=&-\Delta\beta(y) ,\vsp
D(A)&=&\{y\in H^{-1}(\calo)\cap L^1(\calo),\ 
\beta(y)\in H^1_0(\calo)\}\mbox{ and }\vsp
  A&=&\pp\vf,\mbox{ where $\vf(y)=\dd\int_\calo j(y(\xi))d\xi.$}\earr$$
By \eqref{e2.11}, we have also
$$\psi^*(v)=\int_\calo j^*(\beta(\theta))d\xi
+\frac\delta2\,|\theta|^2_{-1},\ v\in L^2(\calo),$$where $\theta\in H^{-1}(\calo)\cap L^1(\calo),$
$$\barr{rcll}
2\delta \theta-\Delta\beta(\theta)&=&v&\mbox{ in }\calo,\vspp
\theta&=&0&\mbox{ on }\pp\calo,\earr$$and $|\cdot|_{-1}$ is the norm of $H^{-1}(\calo)$. Then we have
$$\barr{ll}
\wt G_1(y,y_1)=\!\!\!&\dd\E\int^T_0 \(\int_\calo j(e^{W(t)}(y(t)+x))d\xi+
\dd\frac\delta2\,|e^{W(t)}(y(t)
+x)|^2_{-1}\)d\xi\,dt\vspp
&+\,\dd\E\int^T_0 \int_\calo e^{W(t)}((\nu+\delta)
(y(t)+x)+\mu  x)e^{W(t)}x\,d\xi\,dt\vspp
&-\E\dd\int^T_0(e^{W(t)}y(t),e^{W(t)}(1+\mu)x)dt\vspp
&\dd+\,\E\int^T_0 \eta(e^{W}y)dt
+\dd\frac12\,\E|e^{W(T)}y(T)|^2_{-1}\vsp
&+(e^{W(T)}y_1,e^{W(T)}x)_{-1},\earr$$
while
$$G_2(u)=
 \dd\E\int^T_0 \Bigg(\int_\calo j^*(\beta(\widetilde z))d\xi+\frac\delta2\,|\widetilde z(t)|^2_{-1}\Bigg)dt ,$$where
\begin{equation}\label{e3.27}
\barr{rcll}
 \delta\widetilde z-\Delta\beta(\widetilde z)&=&u&\mbox{ in }\calo,\vspp
\widetilde z&=&0&\mbox{ on }\pp\calo.\earr
\end{equation}
(Here, $(\cdot,\cdot)_{-1}$ is the scalar product of $H^{-1}(\calo).$)

Taking into account that $\frac{j^*(r)}{|r|}\to+\9$ as $|r|\to\9$, it follows, as in the previous case, for each $M>0$, the set
$$\left\{\beta(\widetilde z);\ \E\int^T_0\int_\calo j^*(\beta(\widetilde z))dt\,d\xi\le M\right\}$$is weakly compact in $L^1((0,T)\times\calo\times\ooo)$, we infer that
\begin{equation}\label{e3.28}
\ov G(y,y_1,u)=\wt G_1(y,y_1)+ \E\dd\int^T_0 \bigg(\int_\calo j^*(\beta(\widetilde z))d\xi
+\frac\delta2\,|\widetilde z(t)|^2_{-1}\bigg)dt,
 \end{equation}
where $\widetilde z$ is the solution to \eqref{e3.27}. This implies that
$$D(\ov G)\subset\{(y,y_1,u)\in\calh\times L^2(\ooo;H)\times\calx;\,u\in L^1((0,T)\times\ooo;\calz)\}.$$
Here $\calz=(-\Delta)\1(L^1(\calo))\subset W^{1,p}_0(\calo),\ 1\le p<\frac d{d-1}$, where $\Delta$ is the Laplace operator with homogeneous Dirichlet conditions and $$D(\ov G)=\{(y,y_1,u);\ \ov G(y,y_1,u)<\9\}.$$

We define, as above, the solution to \eqref{e3.25} as $X^*=e^Wy^*$, where $(y^*,y^*_1,u^*)$ is the solution to the minimization problem
\begin{equation}\label{e3.29}
\barr{r}
{\rm Min} \bigg\{\ov G(y,y_1,u); \dd \frac{dy}{dt}{+}
(\mu{+}\nu{+}\delta)(y{+}x){+}e^{-W}u{=}0,\, y(0)=0,\vsp y(T){=}y_1,\
(y,y_1,u)\in\calh\times L^2(\ooo;H)\times\calx
\bigg\}\!\!.\earr
\end{equation}
(Here, $\frac{dy}{dt}$ is taken in sense of distributions, i.e., in $\cald'(0,T;H)$.) We have, therefore,
\begin{theorem}\label{t3.2a} Assume that $x\in L^2(\calo)$. Then equation \eqref{e3.25} has a unique variational solution $X^*$,
$$\barr{c}
X^*\in L^2((0,T)\times\calo\times\ooo);\,\vf(X^*)
\in L^1((0,T)\times\calo\times\ooo),\vsp
e^{-W}X\in W^{1,1}([0,T];W^{1,1}_0(\calo)),\ \pas\earr$$Moreover, the process $t\to X^*(t)$ is pathwise $W^{1,1}_0(\calo)$-valued continuous on~$(0,T).$
\end{theorem}

\begin{remark}\rm A different treatment of equation \eqref{e3.25} under the general assumptions \eqref{e3.26} was developed in \cite{4a} (see also \cite{4aa}, Ch. 5).\end{remark}

\section{Stochastic variational inequalities}
\setcounter{equation}{0}

Consider the stochastic differential equation
\begin{equation}\label{e4.1}
\barr{l}
dX+A_0X\,dt+N_K(X)dt+\lbb X\,dt\ni X\,dW,\ t\in(0, T),\vspp
X(0)=x,\earr\end{equation}in a real Hilbert space $H$ with the scalar product $(\cdot,\cdot)$ and the norm $|\cdot|$. Assume that $x\in H$ and
\begin{itemize}
\item[(j)] {\it $A_0:D(A_0)\subset H\to H$ is a linear self-adjoint, positive definite operator in $H$.}
\item[(jj)] {\it $W$ is the Wiener process \eqref{e1.2} and $\lbb>\nu$.}
\item[(jjj)] {\it $K$ is a closed, convex subset of $H$ such that $0\in K$, $(I+\lbb A_0)^{-1}K\subset K$, $\ff \lbb>0.$ } \end{itemize}
    Here, $N_K:H\to 2^H$ is the normal cone to $K$, that is,
\begin{equation}\label{e4.2}
N_K(u)=\{\eta\in H;\ (\eta,u-v)\ge0,\ \ff v\in K\}.\end{equation}

By the transformation \eqref{e2.1}, equation \eqref{e4.1} reduces to the nonlinear random differential equation
\begin{equation}\label{e4.3}
\barr{l}
\dd\frac{dy}{dt}+e^{-W}A_0(e^W(y+x))+e^{-W}N_K(e^W(y+x))+\mu(y+x)=0,\\\hfill t\in(0,T),\vspp
y(0)=0.\earr\end{equation}
(We note that, if $W(t)=\sum\limits^N_{j=1}\mu_j\b_j(t)$, then \eqref{e4.3} reduces to a deterministic variational inequality.)

To represent this problem as an optimization problem of the form (P), we~set
$$\vf(u)=\frac12\ (A_0u,u)+I_K(u),\ \ff u\in H,$$where $I_K$ is the indicator  function $$I_K(u)=\left\{\barr{ll}
0&\mbox{ if }u\in K,\\
+\9&\mbox{ otherwise}.\earr\right.$$
The function $\vf:H\to]-\9,+\9]$ is convex  and lower semicontinuous. Then, by \eqref{e2.6}, \eqref{e2.17}, \eqref{e2.18}, we have
\begin{eqnarray}
\wt G_1(y,y_1)&=&\E\int^T_0 \bigg(\frac12\,(A_0(e^{W(t)}(y(t)+x)),e^{W(t)}(y(t)+x))\ \ \ \label{e4.4}\\
&&+\,e^{2W}((\nu+\delta)(y+x)+\mu y)x
+\dd\frac\delta2\,|e^{W(t)}(y(t)+x)|^2_H\nonumber\\ &&+I_K(e^{W(t)}(y(t)+x))
+\eta(e^{W(t)}y(t))\bigg)dt\nonumber\\
&&-\E\dd\int^T_0(e^Wy,e^W(1+\mu)y)dt\vspp
&&+\dd\frac12\,\E|e^{W(T)}y_1|^2_H
+\E(e^{W(T)}y_1,e^{W(T)}x),\nonumber\\
G_2(u)&=&\E\int^T_0  \psi^*(u(t) )dt
 , \label{e4.5}
\end{eqnarray}where, by \eqref{e2.11}-\eqref{e2.11a}, we have
\newpage
\begin{equation}\label{e4.5a}
\barr{ll}
\psi^*(e^Wu)\!\!\!
&=\dd\sup\Bigg\{(e^Wu,v)-\frac12\,(A_0v,v)-
\frac\delta2\,|v|^2;\,v\in K\Bigg\}\\
&=(e^Wu,z)-\dd\frac12\,(A_0z,z)-\dd\frac\delta2|z|^2,\earr,
\end{equation}where $A_0z+\delta z+N_K(z)\ni e^Wu.$ (We note that, by (iii), $z$ is uniquely defined.)

By \eqref{e4.4}-\eqref{e4.5}, we see that
$$\ov G(y,y_1,u)=\wt G_1(y,y_1)+\frac12\,\liminf_{n\to\9}
\E\int^T_0\((A_0z_n,z_n)+\frac\delta2\,|z_n|^2\)dt,$$
where
\begin{equation}\label{e4.6-a}
\barr{l}
A_0z_n+\delta z_n+N_K(z_n)\ni  u_n,\ e^W By_n+u_n=0,\vspp
 y_n\to y\mbox{ in }\calh,\ y_n(T)\to y_1\mbox{ in }L^2(\ooo;H),\ u_n\to u\in\calx.\earr\end{equation}This  yields
\begin{equation}\label{e4.7-b}
\E\int^T_0|A^{\frac12}_0z_n|^2 dt\le C<\9,\ \ff n\in\nn,\end{equation}
and, therefore, we have
\begin{equation}\label{e4.9}
\barr{rl}
\ov G(y,y_1,u)\!\!\!
&=\wt G_1(y,y_1)+\dd\frac12\,\E\int^T_0
(|A^{\frac12}_0z|^2+\delta|z|^2)dt,\vspp
z\!\!\!&=w-\dd\lim_{n\to\9}z_n\mbox{ in }L^2((0,T)\times\ooo;V),\earr
\end{equation}where $V=D(A^{\frac12}_0).$
We note that $D(\wt G_1)\subset L^2((0,T)\times\ooo;V).$

We may conclude, therefore, by Theorem \ref{t2.4-b} that

\begin{theorem}\label{t4.1-b} Under hypotheses {\rm(j)--(jjj)}, there is a unique variational solution $X^*(t)\in K$, a.e. $t\in(0,T)$, $X^*\in L^2((0,T)\times V;\ooo)$ to equation \eqref{e4.1}.
\end{theorem}

More insight into the problem can be gained in the following two special cases.

\paragraph{Stochastic parabolic variational inequalities} \ \sk

\n The stochastic differential equation
\begin{equation}\label{e4.18}
\barr{l}
dX-\Delta X\,dt+\lbb X\,dt+N_K(X)dt\ni X\,dW\mbox{ in }(0,T)\times\calo,\vsp
X(0)=x\mbox{ in }\calo,\vsp
X=0\mbox{ on }(0,T)\times\pp\calo,\earr\end{equation}
 where $N_K(X){\subset} L^2(\calo)$ is the normal cone to the closed convex set $K$ of $L^2(\calo)$,
 \newpage
$$K=\{z\in L^2(\calo);\ z\ge0,\mbox{ a.e. in }\calo\},\ \alpha\in\rr,$$can be   treated following the above infinite-dimensional scheme in the space $H=L^2(\calo)$, where $A_0u=-\Delta u,$ $u\in D(A_0)=H^1_0(\calo)\cap H^2(\calo)$.

Then the variational solution to \eqref{e4.18} is defined  by $X=e^Wy$, where $y$ is given by \eqref{e4.8-b} and $\ov G$ is given by
$$\ov G(y,y_1,u)=\wt G_1(y,y_1)+\frac 12\,\E\int^T_0\int_\calo(|\nabla z^2|^2+\delta|z|^2)d\xi\,dt,$$
where $G_1$ is defined y \eqref{e4.4} and $z=w-\lim\limits_{n\to\9}z_n$ in $L^2((0,T)\times\ooo;H^1_0(\calo))$,
\begin{equation}\label{e4.11-b}
\barr{c}
-\D z_n+\delta z_n+\eta_n=u_n,\ e^W \calb y_n+y_n=0,\vsp
\eta_n\in N_K(z_n),\ \E\dd\int^T_0\int_\calo|\nabla z_n|^2d\xi\,dt\le C,\ \ff n.\earr
\end{equation}
Since $u_n\to  u$ in $\cald'(0,T;L^2(\calo))$ and $\eta_n(t,\xi)\le0$ a.e. $(t,\xi)\in(0,T)\times\calo,$  by \eqref{e4.11-b}, we infer that
$$-\Delta z+\delta z+\eta= u\mbox{ in }\cald'((0,T)\times\calo),$$
where  $\eta,u $ are in $\mathcal{M}((0,T)\times\calo)$  the space of bounded measures on \mbox{$(0,T)\times\calo$.} If we denote by $\eta_a,u_a\in L^1((0,T)\times\calo)$ the  absolutely con\-ti\-nuous parts of $\eta$ and $u$, we get
 $$\barr{c}
 -\Delta z+\delta z+\eta_a= u_a\mbox{ in }L^1(\calo),\vsp
 z\in H^1_0(\calo)\mbox{ and }\eta_a(t,\xi)=0,\mbox{ a.e.  on }[z(t,\xi)>0]\vsp
 \eta_a(t,\xi)\ge0\mbox{, a.e. on }[z(t,\xi)=0].\earr$$
Then the process $X=e^W(y+x)$ is the variational solution to \eqref{e4.18} and so, by Theorem \ref{t4.1-b}, we have

\begin{corollary}\label{c4.2}There is a unique variational solution $X{\in} L^2((0,T){\times}\ooo;H^1_0(\calo)),$ $X\ge0$, a.e. on $(0,T)\times\ooo.$
\end{corollary}

 \paragraph{Finite   dimensional   stochastic variational inequalities}\ \sk

   \n Consider equation \eqref{e4.1} in the special case $K\subset\rr^d$, ${\rm int}\ K\ne\emptyset,$ $0\in K$, $W=\sum\limits^N_{i=1}\mu_i\beta_i $ and $A_0\in L(\rr^d,\rr^d)$, $A_0=A^*_0$. Then, as easily seen by \eqref{e2.12}, we have
   \newpage
\begin{equation}\label{e4.9a}
\psi^*(u)\ge\alpha_1|u|-\alpha_2,\ \ff u\in\rr^d.\end{equation}
Let $z_n$ be the solution to (see \eqref{e4.6-a})

 \begin{equation}\label{e4.14}
A_0z_n+\delta z_n+N_K(z_n)\ni u_n.\end{equation}
Since, by \eqref{e4.9a}-\eqref{e4.14}, the sequence $\{ u_n\}$ is bounded in $L^1 ((0,T)\times\ooo,\rr^d)$, it follows that it is weak-star compact in $\mathcal{M}(0,T;\rr^d)$, $\ff\vp>0$, and so $u\in\mathcal{M}(0,T;\rr^d)$. (Here, $\mathcal{M}(0,T;\rr^d)$ is the space of $\rr^d$-valued bounded measures on $(0,T)$. Letting $n\to\9$ in \eqref{e4.14}, we get
\begin{equation}\label{e4.11-a}
A_0z+\delta z+\zeta= u,\end{equation}
where $ u\in \mathcal{M}(0,T;\rr^d)$, $\ff \vp>0,$ and $\zeta\in\mathcal{M}((0,T);\rr^d)$,  $\pas$ By the Lebesgue decomposition theorem, we have
$$\barr{ll}
\zeta=\zeta_a+\zeta_s,&\zeta_a\in L^1 (0,T;\rr^d),\vspp
u=u_a+u_s,&u_a\in L^1 (0,T;\rr^d),\earr$$where $u_s$ and $\zeta_s$ are singular measures and $\zeta_a\in N_K(z).$  Hence, by \eqref{e4.11-a}, we  have
\begin{equation}\label{e4.15prim}
z=(A_0+\delta I+N_K)^{-1}( u_a)=F(u_a),\ \zeta_s= u_s.\end{equation}
As a matter of fact, the singular measure $\zeta_s$ belongs to the normal cone $N_{\mathcal{K}}(z)\subset\mathcal{M}(0,T;\rr^d)$ to the set ${\mathcal{K}}=\{\wt z\in C([0,T];\rr^d);$ $\wt z(t)\in K,$ \mbox{$\ff t\in[0,T]\}$} and it is concentrated on the set of $t$-values for which $z(t)$  defined by \eqref{e4.15prim} lies on the boundary $\pp K$ of $K$.

By \eqref{e2.20a}-\eqref{e2.21}, we have
\begin{equation}\label{e4.15}
\ov G(y,y_1,u)
 =\wt G_1(y,y_1)+\E\dd\int^T_0
 \Bigg(\frac12(A_0F(u_a),F(u_a))
 +\frac\delta2\,|F(u_a)|^2\Bigg)dt,
 \end{equation}where $\wt G_1$ is given by \eqref{e4.4} and $y\in\calh$ is solution to the equation

\begin{equation}\label{e4.16}
\barr{ll}
y=y_a+y_s,\ y_a\in AC([0,T];\rr^d),\ y_s\in BV([0,T];\rr^d),\ \pas,\vsp
\dd\frac{dy_a}{dt}+(\mu+\nu+\delta)
(y_a+x)+e^{-W}u_a
=0,\mbox{ a.e. on }(0,T),\vspp
y_a(0)=0,\vsp
\dd\frac{dy_s}{dt}+ e^{-W}u_s=0\mbox{ in }\cald'(0,T;\rr^d),\earr\end{equation}
where $BV([0,T];\rr^d)$ is the space of functions with  founded variations on $[0,T]$. We note that,  by \eqref{e4.15}, it follows also that
$$\barr{r}
D(\ov G) \subset \{(y,y_1,u) \in \calh\times L^2(\ooo;H)
 \times \calx;\ y \in  BV([0,T];\rr^d),\  \pas,\vsp F(u_a) \in  L^2((0,T) \times \ooo \times \rr^d)\},\earr$$where $D(\ov G)=\{(y,y_1,u);\ G(y,y_1,u)<\9\}.$
We have, therefore,

\begin{theorem}\label{t4.1} The minimization problem
\begin{equation}\label{e4.17}
{\rm Min}\{\ov G(y,y_1,u) ;\,(y,y_1,u)\in\calh\times L^2(\ooo;\rr^d)\times\calx,\ \mbox{subject to \eqref{e4.16}}\}\end{equation}
has a unique solution $(y^*,y^*_1)\in \calh\times L^2(\ooo;\rr^d)$ satisfying \eqref{e4.16}. The process $X^*=e^Wy^*$ is  the solution to the variational solution to \eqref{e4.15}.\end{theorem}

\begin{remark}\rm Since $y^*\in BV([0,T];\rr^d)$  and, as seen by \eqref{e4.16}, the singular  measure $\zeta_s=u_s\ne0$, it follows that the process $X^*$ is not pathwise continuous on $[0,T]$. However, by the Lebesgue decomposition, we have, $\pas$,
 $X^*(t)=X^*_a(t)+X^*_1(t)+X^*_2(t),\ \ff t\in[0,T],$
where $t\to X^*_a(t)e^{-W(t)}$ is absolutely continuous, $X^*_1$ is a jump function and $X^*_2$ is a singular function, that is, $X^*_2=e^Wy_2$, where $\frac{dy_2}{dt}=0$. a.e.\end{remark}

\n{\bf Acknowledgement.} This work was supported by the DFG through CRC~1283. V. Barbu was also partially supported by CNCS-UEFISCDI (Romania) through the project PN-III-P4-ID-PCE-2016-0011.


\begin{thebibliography}{nn}

\bibitem{1}  Barbu, V., {\it Nonlinear Differential Equations of Monotone Type in Banach Spaces}, Springer Monographs in Mathematics, Springer, New York, 2010.\vspace*{-2mm}
\bibitem{2}  Barbu, V., A variational approach to stochastic nonlinear problems, {\it J.~Math. Anal. Appl.}, 384 (2011), 2-15.\vspace*{-2mm}
\bibitem{3}  Barbu, V., Optimal control approach to nonlinear diffusion equations driven by Wiener noise, {\it  J. Optim. Theory Appl.}, 153 (2012), 1-26.\vspace*{-2mm}
\bibitem{3-a}  Barbu, V., A variational approach to   nonlinear stochastic diffrential equations with linear multiplicative noise (submitted).\vspace*{-2mm}
\bibitem{4b} Barbu, Existence for nonlinear finite dimensional stochastic differential equations of subgradient type, {\it Mathematical Control and Related Fields} (to appear).\vspace*{-2mm}
\bibitem{4} Barbu, V., Brzezniak, Z., Hausenblas, E., Tubaro, L., Existence and convergence results for infinite dimensional nonlinear stochastic equations with multiplicative noise, {\it Stoch. Processes and Their Appl.}, 123 (2013), 984-951.\vspace*{-2mm}

\bibitem{4a} Barbu, V., Da Prato, G., R\"ockner, M., Existence of strong solutions for stochastic porous media equations under general monotonicity conditions, {\it Ann. Probab.}, 37 (2) (2009),  428-452.\vspace*{-2mm}

    \bibitem{4aa} Barbu, V., Da Prato, G., R\"ockner, M., {\it Stochastic Porous Media Equations}, Lecture Notes in Mathematics, 2163, Springer, 2016.\vspace*{-2mm}


\bibitem{5} Barbu, V., R\"ockner, M., Stochastic variational inequalities and applications to the total variation flow perturbed by linear multiplicative noise, {\it Archive Rational Mech. Anal.}, 209 (2013), 797-834.\vspace*{-2mm}
\bibitem{5a}  Barbu, V., R\"ockner, M., An operatorial approach to  stochastic partial differential equations driven by linear multiplicative noise, {\it J. European Math. Soc.}, 17 (2015), 1789-1815.\vspace*{-2mm}
\bibitem{6} Brezis, H., Ekeland, I., Un principe variationnel associ\'e \`a certains \'equations paraboliques, le cas ind\'ependent du temps, {\it C.R. Acad. Sci. Paris}, 282 (1976), 971-974.\vspace*{-2mm}
    \bibitem{11-a} Brooks, J.K., Dinculeanu, N., Weak compactness in spaces of Bochner integrable functions and applications, {\it Advances in math.}, 24 (1977), 172-188.\vspace*{-2mm}
\bibitem{7}   Da Prato, G., Zabczyk, J., {\it Stochastic Equations in Infinite Dimensions}, 1992. Second  Edition, Cambridge University Press, Cambridge, 2008.\vspace*{-2mm}
\bibitem{10-a} Gess, B., R\"ockner, M., Stochastic variational inequalities and regularity for degenerate stochastic partial differential equations, {\it Trans. Amer. Math. Soc.}, 369 (2017), 3017-3045.\vspace*{-2mm}
\bibitem{8} Krylov, N.V., Rozovskii, B.L., Stochastic evolution equations, {\it J. Soviet Math.}, 16 (1981), 1233-1277.\vspace*{-2mm}
\bibitem{13-b} Liu, W., R\"ockner, M., {\it Stochastic Partial Differential Equations: An Introduction}, Springer, 2015.\vspace*{-2mm}
\bibitem{9} Pardoux, E., Equations aux d\'eriv\'ees partielles stochastiques nonlin\'eaires monotones, Th\`ese, Orsay, 1972.\vspace*{-2mm}
\bibitem{10} Prevot, C., R\"ockner, M., {\it A Concise Course on Stochastic Partial Dif\-fe\-ren\-tial Equations}, Lecture Notes in Mathematics, 1905, Springer, Berlin 2007.\vspace*{-2mm}
    \bibitem{16-b} Rockafellar, R.T., Integrals which are convex functionals,   {\it Pacific J. Math.}, 24  (1968), 525-539.\vspace*{-2mm}
\bibitem{11} Visintin, A., Extension of the Brezis-Ekeland-Nayroles principle to monotone operators, {\it Adv. Math. Sci. Appl.}, 18 (2008), 633-680.

\end{thebibliography}
\end{document}